\documentclass[12pt,twoside]{extarticle}
\usepackage{mathtext}
\usepackage{titlesec}
\usepackage[cp1251]{inputenc}
\usepackage[T2A]{fontenc}
\usepackage[russian,english]{babel}
\usepackage{fancyhdr}
\usepackage{amsthm}
\usepackage{dsfont}
\usepackage{indentfirst}
\usepackage{mathrsfs}
\usepackage{array}
\usepackage{amssymb}
\usepackage{amsmath}
\usepackage{eucal}
\usepackage{dubloper}
\usepackage{bm}
\usepackage{epsfig}
\usepackage{wrapfig}

\titleformat{\section}[block]{\large\bfseries\filcenter}{\large\bfseries\thesection. }{0pt}{}
\titleformat{\subsection}[block]{\bfseries\filcenter}{\bfseries\thesubsection. }{0pt}{}
\titleformat{\subsubsection}[block]{\bfseries\filcenter}{\bfseries\thesubsubsection. }{0pt}{}

\newcounter{mytheorem}[section]

\edef\lim{\lim\limits}
\renewcommand{\frac}[2]{\dfrac{#1}{#2}}
\numberwithin{equation}{section}
\numberwithin{figure}{section}
\begin{document}

\setlength{\abovedisplayskip}{2mm} 
\setlength{\belowdisplayskip}{2mm} 
\parindent=8,8mm
\renewcommand{\headrulewidth}{0pt}
\edef\sum{\sum\limits}

\thispagestyle{empty}

\begin{center}

\

\

\

\

\textbf{\Large On permutizers of subgroups}

\textbf{\Large of finite groups}

\bigskip

\textbf{\bf A.\,F.\,Vasil'ev, V.\,A.\,Vasil'ev, T.\,I.\,Vasil'eva}

\bigskip

\textbf{Dedicated to Professor Victor D. Mazurov}

\textbf{in honor of his 70th birthday}

\end{center}

\bigskip

\begin{center}
{\bf Abstract}
\end{center}

\medskip

\noindent\begin{tabular}{m{8mm}m{138mm}}

&{\small~~~
Let $H$ be a subgroup of a group $G$. The permutizer of $H$ in $G$ is the
subgroup $P_G(H)=\langle x\in G\ |\ \langle x\rangle H=H\langle x\rangle
\rangle$.
A subgroup $H$ of a group $G$ is called permuteral in $G$, \ if
\ $P_G(H)=G$, \ and \ strongly \ permuteral \ in \ $G$, \ if
\ $P_U(H)=U$ \ whenever $H\leq U\leq G$.
Finite groups with given systems of permuteral and strongly permuteral subgroups
are studied. 
New characterizations of $\rm w$-supersoluble and supersoluble
groups are received.

{\bf Keywords:} 
permutizer of a subgroup, supersoluble group,
$\rm w$-supersoluble group, $\Bbb{P}$-subnormal subgroup, 
\ permuteral subgroup, \ strongly permuteral subgroup, \ product of
subgroups}

\smallskip

MSC2010\  20D20, 20D40 


\end{tabular}

\bigskip

{\parindent=0mm
\textbf{\bf {\large Introduction}}
}
\bigskip

All groups considered in this paper are finite. 
In the theory of groups the normalizer of a subgroup is a classical concept, 
about which there are many well-known results. 
For example, 

\medskip

{\it The following statements about a group $G$
are equivalent$:$

$(1)$ $G$ is nilpotent$;$

$(2)$ $H < P_G(H)$ for every $H < G$ $($the normalizer  condition$);$

$(3)$ $N_G(M)=G$ for every maximal subgroup $M$ of $G$ $($the maximal 
normalizer  condition$);$

$(4)$ $N_G(P)=G$ for every Sylow subgroup $P$ of $G;$

$(5)$ $N_G(S)=G$ for every Hall subgroup $S$ of $G;$

$(6)$ $G=AB$ where $A$ and $B$ are nilpotent subgroups of $G$ and
$N_G(A)=N_G(B)=G$. }

\medskip

A natural generalization of subgroup's normalizer is the concept of the 
permutizer of a subgroup introduced  in 
\cite{Wei}.

\medskip

{\bf Definition 1 \cite[p.~27]{Wei}.} {\it Let $H\leq G$. The permutizer of $H$ in $G$
is the subgroup 
$
P_G(H) = \langle x\in G\ |\ \langle x\rangle H=H\langle x\rangle \rangle.
$
}
\medskip

Replacing in (2)--(6) the normalizer of a subgroup by its  permutizer, we obtain
the following interesting problems.
  
\medskip

{\bf Problem 1.} {\it Describe all groups $G$, satisfying the 
permutizer condition, i.e. $G:$ \ \ \ \ \ \ \ $H <P_G(H)$ for every $H < G$.} 

\medskip
                   
This problem was investigated by W.E.\,Deskins and P.\,Venzke \cite[pp.~27--29]{Wei}, 
J.\,Zhang \cite {Zha}, J.C.\,Beidleman and D.J.S.\,Robinson \cite {BR},
A.\,Ballester-Bolinches and R.\,Esteban-Romero \cite {BE} and others.

\medskip

{\bf Problem 2.} {\it Describe all groups $G$, satisfying the 
maximal permutizer condition, i.e. $G:$ $N_G(M)=G$ for every maximal subgroup $M$ 
of $G$.}

\medskip

This problem was considered by W.E.\,Deskins and P.\,Venzke \cite[pp.~27--29]{Wei},  X.\,Liu and
Ya.\,Wang \cite{Xia}, Sh.\,Qiao, G.\,Qian and Ya.\,Wang 
\cite{Shouhong} and others.

\medskip
In order to briefly formulate assertions (4)--(6) in terms of the permutizers of subgroups,
 we need to introduce the following definition.

\medskip

{\bf Definition 2.} {\it Let $H$ be a subgroup of a group $G$. We say that

$(1)$ $H$ is permuteral in $G$, if $P_G(H)=G;$

$(2)$ $H$ is strongly permuteral in $G$, if $H\leq U\leq G$ then
$P_U(H) = U$.}

\medskip

There exists groups which have permuteral but not strongly permuteral subgroups.
For example, it's easy to check that in the group $G=PSL(2, 7)$ a Sylow
3-subgroup $Z_3$ is permuteral in $G$. Since $Z_3\leq U\leq G$, where
$U$ is isomorphic to the alternating group $A_4$ of degree 4 and  $P_U(Z_3)=Z_3$, then $Z_3$
is not strongly permuteral in $G$.

\medskip
We note the following problems. 

\medskip

{\bf Problem 3.} {\it Describe all groups $G$ such that$:$

$a)$ every Sylow $($Hall$)$  subgroup of $G$ is  permuteral in $G;$

$b)$ every Sylow $($Hall$)$ subgroup of $G$ is  strongly permuteral in 
$G$.}

\medskip

{\bf Problem 4.} {\it Describe all groups $G=AB$  where   $A$ and $B$ are permuteral 
$($strongly permuteral$)$ nilpotent subgroups of $G$.}

\medskip

This paper is devoted the solution of the problems 3 and 4.



\bigskip

{\parindent=0mm
\textbf{\bf {\large 1. Preliminary results}}
}

\bigskip

We use the notation and terminology from \cite{Shem, DH}.
We recall some concepts significant in the paper.

Let $G$ be a group. If $H$ is a subgroup of $G$, we write $H\leq G$ and if $H\not= G$,
we write $H < G$.
We denote by $|G|$ the order of $G$;
by $\pi(G)$ the set of all distinct prime divisors of the order of $G$;
by ${\rm {Syl}}_p(G)$ the set of all Sylow $p$-subgroups of $G$;
by ${\rm {Syl}}(G)$ the set of all Sylow subgroups of $G$;
by ${\rm Core}_G(M)$ the core of subgroup $M$ in $G$, i.e. intersection of all
subgroups conjugated with $M$ in $G$;
by $F(G)$ the Fitting subgroup of $G$;  
by $F_p(G)$ the $p$-nilpotent radical of $G$, i.e.
the product of all normal $p$-nilpotent subgroups of $G$.
$Z_n$ denotes cyclic group of order $n$;
$\Bbb{P}$ denotes the set of all primes; $\pi$ denotes a set of primes;
$\pi'=\Bbb{P}\setminus\pi$;
$\frak{S}$ denotes the class of all soluble groups;
$\frak{U}$ denotes the class of all supersoluble groups;
$\frak{N}$ denotes the class of all nilpotent groups;
$\frak{A}(p-1)$ denotes the class of all abelian groups of exponent dividing $p-1$.

A group $G$ of order $p_1^{\alpha_1}p_2^{\alpha_2}\cdots p_n^{\alpha_n}$ (where
$p_i$ is a prime, $i=1, 2,\ldots, n$) is called {\it Ore dispersive} \cite[p.~251]{Shem},
whenever $p_1>p_2>\cdots >p_n$ and
$G$ has a normal subgroup of order $p_1^{\alpha_1}p_2^{\alpha_2}\cdots p_i^{\alpha_i}$
for every $i=1, 2,\ldots, n$.
{\it A Carter subgroup} of a group $G$ is called a nilpotent subgroup $H$ that
$N_G(H)=H$. A group $G$ is {\it $p$-closed} if $G$ has a normal Sylow $p$-subgroup.

A class of groups $\frak F$ is called a {\it formation} if the following conditions hold:
(a) every quotient group of a group lying in $\frak F$ also lies in
$\frak F$; (b) if $H/A\in \frak F$ and $H/B\in \frak F$ then
$H/A\cap B\in \frak F$.
A formation $\frak{F}$ is called
{\it hereditary} whenever $\frak{F}$ together with every group contains
all its subgroups, and {\it saturated}, if $G/\Phi(G)\in\frak{F}$ implies that $G\in\frak{F}$.
Denote by $G^{\frak {F}}$ the $\frak F$-residual of a group $G$,
i.e. the smallest normal subgroup of $G$ with $G/G^{\frak {F}}\in
\frak {F}$.

A function $f:\Bbb{P}\rightarrow \{$formations$\}$ is called a {\it local function}.
A formation $\frak{F}$ is called {\it local},
if there exists a local function $f$ such that $\frak{F}$
coincides with the class of groups $(G | G/F_p(G)\in f(p)$ for every $p\in\pi(G))$.

\medskip

{\bf Lemma~1.1 \cite[Lemma~3.9]{Shem}.} {\it
If $H/K$ is a chief factor of a group $G$ and $p\in\pi(H/K)$ then $G/C_G(H/K)$
doesn't contain nonidentity normal $p$-subgroup and besides $F_p(G)\leq
C_G(H/K)$.}

\medskip

{\bf Lemma~1.2 \cite[Lemma~1.2]{Shem}.} {\it Let $\frak{F}$ be a nonempty formation,
$K$ be a normal subgroup of a group $G$.
Then $(G/K)^{\frak F}= G^{\frak F}K/K.$}

\medskip

{\bf Theorem~1.3 \cite[A, theorem~2.7~(ii)]{DH}.} {\it Let $G$ be a soluble group. Then
$F(G)/\Phi(G)=C_{G/\Phi(G)}(F(G)/\Phi(G))={\rm {Soc}}(G/\Phi(G))$.}

\medskip

{\bf Definition~1.4 \cite{VasVasTju}.} {\it A subgroup $H$ of a group $G$ is called
$\Bbb{P}$-subnormal in $G$ $($denoted by $H\  \Bbb{P}$-$sn\ G)$,
if either $H=G$, or there exists a chain of subgroups
$H = H_0 < H_1 < \cdots < H_{n-1} < H_{n}
= G$
such that
$| H_{i+1} : H_{i}|$ is a prime for every $i = 0, 1,\ldots, n-1$.
}

\medskip

{\bf Lemma~1.5 \cite[Lemmas~3.1 and 3.4]{VasVasTju1}.} {\it Let $H$ be a
subgroup of a group $G$.
Then$:$ 

$(1)$ If $N \trianglelefteq\ G$ and $H\  \mathbb{P}$-$sn\ G$ then
$(H \cap N)\  \mathbb{P}$-$sn\ N$
and $HN/N\  \mathbb{P}$-$sn\ G/N.$

$(2)$ If $N \trianglelefteq\ G$, $N\leq H$ and $H/N\  \mathbb{P}$-$sn\ G/N$ then
$H\  \mathbb{P}$-$sn\ G.$

$(3)$ If $HN_i\  \mathbb{P}$-$sn\ G$, $N_i \trianglelefteq\ G$, $i=1, 2$
then $(HN_1 \cap HN_2)\  \mathbb{P}$-$sn\ G.$

$(4)$ If $H\  \mathbb{P}$-$sn\ K$ and $K\  \mathbb{P}$-$sn\ G$
then $H\  \mathbb{P}$-$sn\ G.$

$(5)$ If $H\  \mathbb{P}$-$sn\ G$ then $H^x\  \mathbb{P}$-$sn\ G$
for every $x \in G.$

$(6)$ If $ G^\mathfrak{U} \leq H$ then $H\  \mathbb{P}$-$sn\ G.$

$(7)$ If $G$ is soluble, $H\  \mathbb{P}$-$sn\ G$ and $K$ is a subgroup of $G$
then $(H \bigcap K)\  \mathbb{P}$-$sn\ K.$

$(8)$ If $G$ is soluble, $H_i\  \mathbb{P}$-$sn\ G,\ i=1,2$ then
$(H_1 \bigcap H_2)\  \mathbb{P}$-$sn\ G$.}

\medskip

A group $G$ is called {\it ${\rm w}$-supersoluble} \cite{VasVasTju}, if every Sylow
subgroup of $G$ is $\Bbb{P}$-subnormal in $G$.
Denote by ${\rm w}\frak U$ the class of all w-supersoluble groups. 
Observe that  $\frak U\subseteq{\rm w}\frak U$.

The following example \cite{VasVasTju} shows that $\frak U\not={\rm 
w}\frak U$.
Take the symmetric group $S$ of degree 3. According to 
\cite[Chapter B, Theorem 10.6]{DH}, there exists a faithful irreducible 
$S$-module $U$ over the field $F_7$ with 7 elements. Consider the semidirect 
product $G = [U]S$. Since $S$ is a nonabelian group, $G$ is not supersoluble. 
The supersolubility of $G/U$ implies that $H_1 = UG_2, H_2 = UG_3$, 
and $H_3 = UG_7 = G_7$ are $\Bbb{P}$-subnormal subgroups of $G$, 
where $G_p $ is a Sylow $p$-subgroup of $G$ for $p\in \{2,3,7\}$. 
Observe that $H_i$ supersoluble subgroup of $G$ for $ = 1,2,3$. Consequently, 
$G_2\  \Bbb{P}$-$sn\ H_1$ and $G_3\  \Bbb{P}$-$sn\ H_2$. This implies that 
$G_p\  \Bbb{P}$-$sn\ G$ for $p\in \{2,3,7\}$, and so $G\in{\rm w}\frak U$. 

\medskip

We present some properties of $\rm{w}$-supersoluble groups.

\medskip

{\bf Proposition~1.6 \cite[Proposition~2.8]{VasVasTju}.} {\it Every
$\rm{w}$-supersoluble group is Ore dispersive.}

\medskip

{\bf Theorem~1.7 \cite[Theorems~2.7, 2.10]{VasVasTju}.} {\it The class $\rm{w}\frak {U}$
is a hereditary saturated formation and it has a local function $f$ such that
$f(p)=(G\in\frak S | {\rm {Syl}}(G)\subseteq
\frak {A}(p-1))$ for every prime $p$.}

\medskip

{\bf Theorem~1.8 \cite[Theorem~2.13]{VasVasTju}.} {\it Every biprimary
subgroup of a $\rm{w}$-supersoluble
group is supersoluble.}

\medskip

{\bf Theorem 1.9 \cite[I, Theorem~1.4]{Wei}.}  {\it
Let  $H/K$ be a chief $p$-factor of a group $G$. $|H/K|=p$ if and only if
${\rm Aut}_G(H/K)$ is abelian group of exponent dividing $p-1$.}

\medskip

A subgroup $H$ of a group $G$ is called: 1) {\it pronormal in $G$}, if the subgroups  
$H$ and $H^x$ are conjugated in their join $\langle H, H^x\rangle$ for all $x\in G$;
2) {\it abnormal in $G$}, if $x\in \langle H, H^x\rangle$ for all
$x\in G$.

\medskip

{\bf Lemma~1.10 \cite [Lemma~17.1] {Shem}.} {\it
If a subgroup $H$ is pronormal in $G$ then $N_G(H)$ is abnormal in $G$.}

\medskip

{\bf Lemma~1.11 \cite [Lemma~17.2] {Shem}.} {\it
Let $H$ be a subgroup of a group $G$. Then the following conditions are
equivalent$:$

$(1)$ $H$ is abnormal in $G;$

$(2)$ If $H\leq U \leq G$ and $H\leq U\cap U^x$ implies that $x\in U;$

$(3)$ $H$ is pronormal in $G$ and $U=N_G(U)$ for $H\leq U\cap U^x;$

$(4)$ $H$ is pronormal in $G$ and $H=N_G(H)$.}

\medskip

{\bf Lemma~1.12 \cite [Lemma~17.5] {Shem}.} {\it
Let $H$ be a subgroup of a group $G$. Then$:$ 

$(1)$ If $H$ is pronormal in $G$ and $H\leq U\leq G$, then $H$ is pronormal in $U.$

$(2)$ Let $N\trianglelefteq\ G$ and $N\leq H$. Then
$H$ is pronormal in $G$ if and only if
$H/N$ is pronormal in $G/N.$

$(3)$  If $N\trianglelefteq\ G$ and $H$ is pronormal in $G$ then
$HN/N$ is pronormal in $G/N$.}

\medskip

{\bf Lemma~1.13 \cite[Lemma~2]{Hei}.} {\it Assume that $G$ is the product 
of two nilpotent subgroup $A$ и $B$, and that $G$ possesses in addition a 
minimal normal subgroup $N$ such that $N=C_G(N)\not=G$. Then

$(1)$ $A\cap B=1.$

$(2)$ $N\leq A\cup B.$

$(3)$ If $N\leq A$ then $A$ is a $p$-group for some prime $p$, and $B$ is a $p'$-group.}

\medskip

The group $G$ satisfies $E_{\pi}$ if $G$ has at least one Hall $\pi$-subgroup;
it satisfies $D_{\pi}$ if $G$ there exists precisely one conjugacy class of Hall
$\pi$-subgroups and if every $\pi$-subgroup of $G$ is contained in a Hall $\pi$--subgroup of $G$.

\medskip

{\bf Lemma~1.14 \cite{Penn}.} {\it Let $G=AB$ be a group satisfying $D_{\pi}$  and
suppose both $A$ and $B$ satisfy $E_{\pi}$. Then there exists Hall
$\pi$-subgroups $A_{\pi}$ and $B_{\pi}$ of $A$ and $B$ respectively such that
$A_{\pi}B_{\pi}=B_{\pi}A_{\pi}$ is a Hall $\pi$-subgroup of $G$.}

\bigskip

{\parindent=0mm
\textbf{\bf {\large 2. Properties of permutizers of subgroups}}
}

\bigskip

{\bf Lemma 2.1.} {\it Let $H$ be a subgroup of a group $G$. Then$:$

$(1)$ $P_U(H)\leq P_G(H)$ for every subgroup $U$ of $G$ such
that $H\leq U.$

$(2)$ $P_G(H)^g=P_G(H^g)$ for every $g\in G.$

$(3)$ $N_G(H)\leq P_G(H).$

$(4)$ If $N\trianglelefteq\ G$ then $P_G(H)N/N\leq P_{G/N}(HN/N).$

$(5)$ If $N\trianglelefteq\ G$ and $N\leq H$ then
$P_{G/N}(H/N) = P_G(H)/N$.}

\medskip

{\it Proof.}
(1) follows from the definition of $P_G(H)$.

(2) Let $g\in G$. Suppose that $P_G(H) = \langle L\rangle$, where
$L=\{x\in G\ |\ \langle x\rangle H=H\langle x\rangle\}$, and
$P_G(H^g) = \langle K\rangle$, where
$K=\{y\in G\ |\ \langle y\rangle H^g=H^g\langle y\rangle\}$. Clearly,
$P_G(H)^g = \langle L^g\rangle$.

Let's take any $z\in L^g$. Then $z=x^g$ for some $x\in L$.
From $\langle x^g\rangle H^g=\langle x\rangle^g H^g=(\langle x\rangle H)^g=
(H\langle x\rangle)^g= H^g\langle x^g\rangle$ we obtain that 
$z\in K$. Hence $L^g\subseteq K$.

We consider any $y\in K$. From $y^{g^{-1}}\in K^{g^{-1}}$ we have
$\langle y^{g^{-1}}\rangle H=\langle y\rangle^{g^{-1}} (H^g)^{g^{-1}}=
(\langle y\rangle H^g)^{g^{-1}}=(H^g\langle y\rangle)^{g^{-1}}=
H\langle y^{g^{-1}}\rangle$. Hence $y^{g^{-1}}\in L$ and $K\subseteq L^g$.

Thus $P_G(H)^g = \langle L^g\rangle= \langle K\rangle =P_G(H^g)$.

Statements (3)--(5) follows from Lemma~2.4 \cite{Shouhong}. \hspace{\stretch{1}}$\square$

\medskip

It is easy to verify the following result.

\medskip

{\bf Lemma 2.2.} {\it Let $H$ be a subgroup of a group $G$ and
$N\trianglelefteq\ G$. Then$:$

$(1)$ If $H$ is permuteral in $G$ then  $HN/N$ is permuteral in $G/N.$

$(2)$ If $H$ is permuteral in $G$ then  $HN$ is permuteral in $G.$

$(3)$ If $N\leq H$ then $H$ is permuteral in $G$ if and only if $H/N$ is permuteral in $G/N.$

$(4)$ If $H$ is strongly permuteral in $G$ then  $HN/N$ is strongly permuteral in $G/N$.}

\medskip

{\bf Lemma 2.3.} {\it Let $G=HQ$ be a group, where $H\in {\rm Syl}_p(G)$,
$p$ is the largest prime divisor of $|G|$, $Q$ is a cyclic subgroup of $G$.
Then $G$ is $p$-closed.}

\medskip

{\it Proof.} Let $G$ be a group of minimal order for which the lemma is false.
Since $G$ is a product of nilpotent subgroups, by the Theorem of
Kegel-Wilandt \cite{Keg}, \cite{Wiel}
$G$ is soluble. Let $N$ be a minimal normal subgroup of $G$. Then
$G/N$ $p$-closed. Since the class of all $p$-closed groups is a saturated
formation, $N$ is an unique minimal normal subgroup of $G$, $\Phi(G)=1$.
Then there exists a maximal subgroup $M$ of $G$ such that
$G=NM$, where $M\cap N=1$, ${\rm Core}_G(M)=1$ and $N=C_G(N)$.
If $N$ is a $p$-group then $HN/N=H/N\in {\rm Syl}_p(G/N)$. Hence
$H\trianglelefteq\ G$, a contradiction.
Let $N$ be a $q$-group, $q\not=p$. In view of Sylow Theorem
$H^g\leq M$ for some $g\in G$ and
$N\leq Q$.
Then $|N|= q$. Hence $M\simeq G/C_G(N)$ can be
embedded in ${\rm Aut}(Z_q)\simeq Z_{q-1}$,
a contradiction with $p > q$.
\hspace{\stretch{1}}$\square$

\medskip

{\bf Lemma 2.4.} {\it Let $H\in {\rm Syl}_p(G)$ and
$p$ be the largest prime divisor of $|G|$. If $H$ is permuteral in $G$ then
$G$ $p$-closed.}

\medskip

{\it Proof.} Let $x$ be an arbitrary element of a group $G$ such that
$\langle x\rangle H=H\langle x\rangle$. Then  $\langle x\rangle H$ is a subgroup of $G$.
By Lemma~2.3 $H\trianglelefteq\langle x\rangle H$. So
$\langle x\rangle\leq N_G(H)$ and $G=P_G(H)\leq N_G(H)$.
 \hspace{\stretch{1}}$\square$

\medskip

{\bf Lemma 2.5.} {\it If every Sylow subgroup of a group $G$ is permuteral
in $G$ then $G$ is Ore dispersive.}

\medskip

{\it Proof.} We prove the lemma by using induction on $|G|$.
We may assume that $|\pi(G)| > 1$. Let
$|G|=p_1^{n_1}p_2^{n_2}\cdots p_k^{n_k}$, where
 $p_1 > p_2 > \cdots > p_k$, $p_i$ are
primes, $i=1, 2,\ldots, k$. For $P_1\in
{\rm Syl}_{p_1}(G)$ by Lemma~2.4 $P_1\trianglelefteq\
G$. Every Sylow $p_i$-subgroup of a quotient group $G/P_1$ has the form
$P_iP_1/P_1$, where $P_i\in {\rm Syl}_{p_i}(G)$, $i=2,\ldots, k$.
In view of Lemma~2.2(1) $P_iP_1/P_1$ is permuteral in $G/P_1$. By
induction $G/P_1$ is Ore dispersive.
Hence $G$ is Ore dispersive.  \hspace{\stretch{1}}$\square$

\medskip

{\bf Lemma 2.6.} {\it If $G$ is a supersoluble group then every pronormal
subgroup of $G$ is strongly permuteral in $G$. }

\medskip

{\it Proof.} In view of the heredity of $\frak U$ 
and Lemma~1.12(1) it suffices to prove that  every pronormal
subgroup of $G\in \frak U$ is permuteral in $G$.

Let $G$ be a supersoluble group of minimal order such that
$P_G(H)\not=G$ for some pronormal subgroup $H$ of $G$.

We suppose that $\Phi=\Phi(G)\not=1$. Then $G/\Phi\in\frak U$,
by Lemma~1.12(3) 
$H\Phi/\Phi$ is pronormal in $G/\Phi$. By the choice of $G$ it follows
$P_{G/\Phi}(H\Phi/\Phi) =G/\Phi$. In view of Lemma~2.1(5) we conclude that
$P_{G}(H\Phi) =G$. Since $P_G(H)\not=G$, there exists an element $x\in G$ such that $x\not
\in P_G(H)$ and $\langle x\rangle H\Phi= H\Phi\langle x\rangle$. Then $R=
\langle x\rangle H\Phi$ is a subgroup of $G$. If
$R\not= G$ then the choice of $G$ implies that $P_R(H) =R$. So
$x\in R=P_R(H)\leq P_G(H)$, a contradiction. 
So $R=\langle x\rangle H\Phi=G=\langle x\rangle H$.
Hence $x\in P_G(H)$, a contradiction.

Thus $\Phi(G)=1$. Group $G\in\frak U$,
so its commutator subgroup $G'\in\frak N$.
The choice of $G$ implies that $N_G(H)\not=G$.
The abnormality of $N_G(H)$ in $G$
implies $G=G'N_G(H)=F(G)N_G(H)$. By Theorem~1.3 
$F(G)=N_1 \cdots N_k$,
where $N_i$ is a minimal normal subgroup of $G$ and $i=1, \ldots, k$.
The supersolvability of $G$ implies $N_i$ is a cyclic subgroup of a prime order. From $N_iH=HN_i$
we get $N_i \leq P_G(H)$ for any $i=1, \ldots, k$.
So $F(G) \leq P_G(H)$. But then $G=F(G)N_G(H)\leq P_G(H)$,
a contradiction.
\hspace{\stretch{1}}$\square$

\medskip

{\bf Corollary 2.6.1. } {\it If $G$ is a supersoluble group then
every Sylow sub\-group of $G$ is strongly permuteral in $G$.}

\medskip

{\bf Corollary 2.6.2. } {\it If $G$ is a supersoluble group then every
Carter subgroup of $G$ is strongly permuteral in $G$. }

\medskip

{\bf Corollary 2.6.3. } {\it If $G$ is a supersoluble group then
every Hall subgroup of $G$ is strongly permuteral in $G$. }

\medskip
The following example shows that there exists 
supersoluble groups with nonpermuteral subgroups.

{\bf Example 2.7.} Let $G=(Z_4\times Z_2)Z_2=\langle a, b |
a^4=b^4=(ab)^2=(a^{-1}b)^2=1\rangle$ be a group of order 16. Then the subgroup
$H=\langle ba\rangle$
is not permuteral in $G$, since $P_G(H)$ is an elementary abelian
2-group of order 8.

\medskip

{\bf Lemma 2.8.} {\it Let $G$ be a soluble group.
If $H$ is a $\Bbb{P}$-subnormal Hall subgroup of $G$
then $H$ is strongly permuteral in $G$.}

{\it Proof.} In view of the heredity of $\frak S$ 
and Lemma~1.5(7) it suffices to prove that  every  $\Bbb{P}$-subnormal 
Hall subgroup of $G\in \frak S$ is permuteral in $G$.

Let $G$ be a soluble group of minimal order
which has a Hall $\pi$-subgroup $H$,
$H\  \Bbb{P}$-$sn\ G$ and $P_G(H)\not= G$.
Let $N$ be a minimal normal subgroup of $G$. Then
$HN/N$ is a Hall $\pi$-subgroup of $G/N$. By Lemma~1.5(1) 
$HN/N\  \Bbb{P}$-$sn\ G/N$. By the choice of $G$, the Hall
$\pi$-subgroup $HN/N$ is permuteral in $G/N$. By Lemma~2.2(3)
$HN$ is permuteral in $G$. 
Therefore $N$ is a $q$-group for some prime $q\not \in \pi$.
Since $H\  \Bbb{P}$-$sn\ G$, it follows that there exists a maximal subgroup 
$M$ in $G$ such that $H\leq  M$ and $|G : M|$ is a prime.
By Lemma~1.5(7) $H\  \Bbb{P}$-$sn\ M$. The choice of $G$
implies that $M=P_M(H)\leq P_G(H)\not= G$. Therefore $M=P_G(H)$.
Since $G=P_G(HN)$, there exists $x\in G$ such that
$\langle x\rangle HN=HN\langle x\rangle$ and $x\not\in M$.
Hence $P_G(H)=M$ implies that $G= \langle x\rangle HN$.
If $N\leq \Phi(G)$ then
$G=
\langle x\rangle H$. Hence $x\in P_G(H)=M$ that contradicts $x\notin M$.

So $N\not\leq \Phi(G)$. Then there exists a maximal subgroup $W$ in $G$ such that $N\not\leq W$
и $G=NW$. Hence $|G : W|$ is a $q$-number and
$H\leq W^g$ for some $g\in G$. Then
$W^g=P_{W^g}(H)\leq P_G(H)=M$ and 
$G=NM$.

Assume that $HN\not=G$. Then by the choice of $G$ we conclude that
$HN=P_{HN}(H)\leq P_G(H)=M$. We get the contradiction
$G=NM\leq M\not=G$.

Hence $HN=G$. Since $N\cap M=1$, it follows that $H=M$.
Then $|N|=q$. In view of $HN=NH$ we get $N\leq
P_G(H)=M$. Hence $G\leq M\not=G$,
a contradiction. 
\hspace{\stretch{1}}$\square$

\medskip

{\bf Corollary 2.8.1. } {\it If $G$ is a $\rm w$-supersoluble group then every
Sylow subgroup of $G$ is strongly permuteral in $G$.}

\bigskip

{\parindent=0mm

\textbf{\bf {\large 3. Characterizations of $\rm w$-supersoluble and supersoluble groups}}

}
\bigskip

{\bf Theorem 3.1.} {\it A group $G$ is $\rm w$-supersoluble if and only if
every Sylow subgroup of $G$ is strongly permuteral in $G$.}

\medskip

{\it Proof.} Let $\frak F =(G\ |\ $ any Sylow subgroup of the group $G$ is strongly
permuteral in $G$).
Сorollary~2.8.1 implies that ${\rm w}\frak U\subseteq \frak F$.

Let $G$ be a group of minimal order in $\frak F\setminus {\rm w}\frak U$.
Since $P_G(H) = G$ for every Sylow subgroup $H$
of $G$, by Lemma~2.5 $G$ is Ore dispersive. Let $N$ be a minimal
normal subgroup of $G$. Take an arbitrary Sylow
$p$-subgroup $R/N$ of $G/N$.
There exists $P \in \text{Syl}_p(G)$ such that $R/N=PN/N$.
Since $G\in \frak F$, by Lemma~2.2(4) it follows $G/N\in \frak F$.
By the choice of $G$ the quotient group $G/N\in {\rm w}\frak U$.
Since ${\rm w}\frak U$ is a saturated formation
by 
Theorem~1.7, so $N$ is an unique minimal normal subgroup of $G$ and $\Phi(G)=1$.
Let $|G|=p_1^{n_1}p_2^{n_2}\cdots p_k^{n_k}$, where $p_i$ are
primes, $i=1, 2,\ldots, k$ and $p_1 > p_2 > \cdots > p_k$. Denote by $P_i$ Sylow
$p_i$-subgroup of $G$, $i=1, 2,\ldots, k$.
Then $P_1\trianglelefteq\ G$ and  $N\leq P_1$. Note that
$P_1\  \Bbb{P}$-$sn\ G$.

Let's denote $H_i=P_iP_1$ for $i\in\{2,\ldots, k\}$.

If $H_i\not= G$ for every $i\in\{2,\ldots,
k\}$, then $H_i\in \frak F$. By the choice of $G$, $H_i$ is $\rm w$-supersoluble.
The heredity of ${\rm w}\frak U$ implies that
$P_iN\in {\rm w}\frak U$. So $P_i\  \Bbb{P}$-$sn\ P_iN$.
Since $G/N\in {\rm w}\frak U$, $P_iN\  \Bbb{P}$-$sn\ G$ by Lemma~1.5(2).
By Lemma~1.5(4) we get $P_i\  \Bbb{P}$-$sn\ G$. Hence $G\in {\rm w}\frak U$.
It contradicts the choice of
$G$.

Hence $H_i=G$ for some $i\in\{2, \ldots, k\}$. Since $\Phi(G)=1$, so
$G=MN$, where $M$ is some maximal subgroup of
$G$, $M\cap N=1$, $N=C_G(N)$. Since $G/N\simeq M$, by Lemma~1.1 it implies
$P_1\cap M=1$. Therefore $N=P_1$ and $M=P_i$.
In view of $P_G(P_i)=G$, there exists an element $y$ of $G$ such that
$y\not\in P_i$, $\langle y\rangle P_i=P_i\langle y\rangle$.
Then $G=\langle y\rangle P_i$. Hence $|N|=p_1$.  
So $G$ is supersoluble,
a contradiction. 
 \hspace{\stretch{1}}$\square$

\medskip


Recall \cite[p. 519]{DH} that the nilpotent length of a soluble group $G$, 
is the smallest integer $l$ such that $F_l(G)= G$, where 
the subgroups $F_i(G)$ are determined recursively as $F_0(G) = 1$ and 
$F_i(G)/F_{i-1}(G) = F(G/F_{i-1}(G))$ for all $i \geq 1$. 
The proposition~2.5 \cite{VasVasTju} shows that 
the nilpotent length of a $\rm w$-supersoluble group cannot be bounded by a fixed 
integer $n$. Since the commutator subgroup 
of every supersoluble group is nilpotent, the nilpotent length of a supersoluble 
group is at most 2, i.e. every supersoluble group is metanilpotent. 
 
\medskip

{\bf Theorem 3.2.} {\it Let $G$ be a metanilpotent group. Then the following statements 
are equivalevt$:$

$(1)$  $G$ is supersoluble$;$

$(2)$  Every Sylow subgroup of $G$ is strongly permuteral in $G;$

$(3)$  Every Sylow subgroup of $G$ is permuteral in $G$.}

\medskip

{\it Proof.} $(1)\Rightarrow (2)$ by Corollary~2.6.1.

$(2)\Rightarrow (3)$ obviously. 

$(3)\Rightarrow (1)$. 
Let $\frak F =(G\ |\ $ the group $G$ is metanilpotent
and every Sylow subgroup of $G$ is permuteral in $G$).

Let $G$ be a group of minimal order in
$\frak F\setminus \frak U$. The nilpotency of
$\frak N$-residual $G^{\frak N}$ and the quotient group
$G/G^{\frak N}$ implies the solvability of $G$. Let $N$ be a minimal normal
subgroup of $G$. By Lemma~1.2 and the nilpotency of $G^{\frak N}$
implies that $(G/N)^{\frak N}=
G^{\frak N}N/N\simeq G^{\frak N}/G^{\frak N} \cap N$ is nilpotent.
For $R/N\in{\rm {Syl}}_p(G/N)$ ($p$ is any prime of $\pi(G)$) there exists
$P\in{\rm {Syl}}_p(G)$ such that $R/N=PN/N$. In view of Lemma~2.2(1),
$R/N$ is permuteral in $G/N$. By the choice of $G$, we get $G/N$ is supersoluble.
From the saturation of the formation $\frak U$ of all supersoluble groups, we
conclude that $N$ is an unique minimal normal subgroup of $G$, $\Phi(G)=1$. Then  $G=NM$,
where $M$ is a maximal subgroup of $G$, $N\cap M=1$, ${\rm Core}_G(M)=1$. Since $N$
is an elementary abelian $p$-group for some $p\in\pi(G)$ and $N=C_G(N)$,
by using of Lemma~1.1, we get $N=G^{\frak N}=F(G)$ and
$G/N\simeq M$ is nilpotent $p'$-group. By Lemma~2.5, $G$ is Ore dispersive. Then for
$|G|=p_1^{n_1}p_2^{n_2}\cdots p_k^{n_k}$ ($p_i$ are primes, $i=1, 2,\ldots, k$ and
$p_1 > p_2 > \cdots > p_k$),
Sylow $p_1$-subgroup $P_1\trianglelefteq\ G$. Hence $p=p_1$ and $N=P_1$.

We will fix $i\in\{2, \ldots, k\}$.
Let $P_i\in{\rm {Syl}}_{p_i}(G)$ and $P_i\leq M$.
Since $P_G(P_i)=G$, there exists $x\in G$ such that
$\langle x\rangle P_i=P_i\langle x\rangle$ and $x\notin M$.

If $\langle x\rangle$ is a $p_1'$-group then $\langle x\rangle P_i$ is also
a $p_1'$-group. The solvability of $G$ implies that $\langle x\rangle P_i
\leq M^g$ for some $g\in G$. Then $\langle x\rangle
^{g^{-1}}P_i^{g^{-1}}\leq M$. Since $P_i^{g^{-1}}$ is a Sylow
$p_i$-subgroup of the nilpotent group $M$, we get 
$P_i^{g^{-1}}=P_i$. Hence $g^{-1}\in N_G(P_i)=M$. Therefore $g\in M$.
Hence $x\in \langle x\rangle P_i\leq M^g=M$. We get a contradiction with
$x\notin M$.

Thus $\langle x\rangle$ is not a $p_1'$-group. Let
$\langle z\rangle\in{\rm {Syl}}_{p_1}(\langle x\rangle)$. Clearly,
$\langle z\rangle\in{\rm {Syl}}_{p_1}(\langle x\rangle P_i)$ and
$\langle z\rangle=P_1\cap\langle x\rangle P_i\trianglelefteq\ \langle x\rangle P_i$.
So $\langle z\rangle P_i$ is a subgroup of $\langle x\rangle P_i$.
Since $\langle z\rangle\leq P_1$ and $P_1$ is an elementary
abelian $p_1$-group, we conclude $|\langle z\rangle|=p_1$.

Let's denote $R_i=\langle z\rangle P_i$. The subgroup $P_i$ is maximal in $R_i$.
Since $N_G(P_i)=M$, it follows $\langle z\rangle\not\leq N_{R_i}(P_i)$.
Therefore $N_{R_i}(P_i)=P_i$. We will show $C_i=C_{R_i}(\langle z\rangle)\cap P_i=1$.
Assume that $C_i\not=1$. Then $\langle z\rangle\leq N_G(C_i)$ and
$P_i\leq N_{R_i}(C_i)\leq N_G(C_i)$. Since $M$ is nilpotent,
we get $P_j\leq N_G(C_i)$ for any $j\in\{2, \ldots, k\}, j\not=i$.
So $M\leq N_G(C_i)$ and $M\not= N_G(C_i)$. Hence $C_i\trianglelefteq G$.
Therefore $1\not=C_i\leq{\rm {Core}}_G(M)=1$,
a contradiction. So $C_i=C_{R_i}(\langle z\rangle)\cap P_i=1$.
Then $P_i\simeq R_i/\langle z\rangle=R_i/C_{R_i}(\langle z\rangle)$
can be embedded in ${\rm {Aut}}(Z_{p_1})\simeq Z_{{p_{1}-1}}$.

Thus $P_i\in\frak{A}(p_1-1)$ for all $i\in\{2, \ldots, k\}$. 
Hence the nilpotency of $M$ implies that $M\in\frak{A}(p_1-1)$. 
Since $M\simeq G/N=G/C_G(N)$, $|N|=p_1$ by Theorem~1.9. So $G$ is supersoluble,
a contradiction.
\hspace{\stretch{1}}$\square$

\medskip

{\bf Theorem 3.3.} {\it Let $G$ be a group. Then the following statements 
are equivalevt$:$

$(1)$  $G$ is supersoluble$;$

$(2)$  Every pronormal subgroup of $G$ is strongly permuteral in $G;$

$(3)$  Every pronormal subgroup of $G$ is permuteral in $G;$

$(4)$ Every Hall subgroup of $G$ is strongly permuteral in $G;$

$(5)$ Every Hall subgroup of $G$ is permuteral in $G$. }

\medskip

{\it Proof.} $(1)\Rightarrow (2)$ by Lemma~2.6.

$(2)\Rightarrow (3)$ obviously. 

$(2)\Rightarrow (4)$. Since every Sylow subgroup of $G$ is pronormal in
$G$, by (2) and Lemma~2.5 $G$ is solvable. Then every Hall subgroup of $G$ 
is pronormal in $G$ and by (2) is strongly permuteral in $G$.

$(4)\Rightarrow (5)$ obviously.

$(3)\Rightarrow (5)$. By (3) and Lemma~2.5 $G$ is soluble. Then every Hall subgroup of $G$ 
is pronormal in $G$ and by (3) is permuteral in $G$.

$(5)\Rightarrow (1)$. 
Let $\frak F =(G\ |\ $ every Hall
subgroup of $G$ is permuteral in $G$).

Let $G$ be a group of minimal order in $\frak F\setminus
\frak U$. Then every Sylow subgroup of $G$ is permuteral in $G$.
By Lemma~2.5, $G$ is Ore dispersive.
Let $N$ be a minimal normal subgroup of $G$. For a set of primes $\pi$ take any
Hall subgroup $\pi$-subgroup
$K/N$ of $G/N$. The solvability of $G$ implies that $K/N=
SN/N$ for some Hall $\pi$-subgroup $S$ of $G$. Since
$G\in \frak F$, by Lemma~2.2(1), we conclude $K/N=SN/N$
is permuteral in $G/N$. By the choice of $G$, it follows that $G/N\in\frak
U$.
Since $\frak U$ is a saturated formation, $N$ is an unique minimal normal subgroup
of $G$, $\Phi(G)=1$.
Then there exists a maximal subgroup $M$ in $G$ such that
$G=MN$, $M\cap N=1$, $N=C_G(N)$. In view of $P\trianglelefteq\ G$ for
$P\in {\rm {Syl}}_p(G)$ ($p$ is the largest prime divisor of $|G|$)
we get $N\leq P$. Since $P\cap M\trianglelefteq\ M$,
by Lemma~1.1 it implies that $P\cap M=1$. Therefore $N=P$ and $M$ is Hall
$\omega$-subgroup for $\omega=\pi(G)\setminus\{p\}$.
Then $M$ is permuteral in $G$. So there exists $x\in G$, $x\not\in M$
such that $\langle x\rangle M=G$. Then a Sylow $p$-subgroup
of $\langle x\rangle$ is a Sylow $p$-subgroup of
$G$. Therefore $|N|=p$. Hence  $M\simeq G/C_G(N)$
is isomorphic embedded in ${\rm Aut}Z_p\simeq Z_{p-1}$. So
$G\in\frak U$, 
a contradiction. 
 \hspace{\stretch{1}}$\square$

\medskip

{\bf Corollary 3.3.1 \cite{KnMo}. } {\it If every Hall subgroup of a group 
$G$ is $\mathbb{P}$-subnormal in $G$
then $G$ is supersoluble. }

\medskip

{\it Proof.}  Since every Sylow subgroup 
of a group $G$ is $\mathbb{P}$-subnormal in $G$, 
$G$ is soluble by Proposition~1.6. By Lemma~2.8 and Theorem~3.2 we conclude that $G$
is supersoluble.

\medskip

{\bf Theorem 3.4.} {\it Let $G$ be a group. Then the following statements 
are equivalevt$:$

$(1)$  $G$ is supersoluble$;$

$(2)$  $G=AB$ is the product of strongly permuteral nilpotent subgroups $A$ and $B$ of $G;$

$(3)$  $G=AB$ is the product of permuteral nilpotent subgroups $A$ and 
$B$ of $G$.}


\medskip

{\it Proof.} $(1)\Rightarrow (2)$.
If $G$ is supersoluble, then $G=F(G)H$, where
$H$ is a Carter subgroup of $G$. Subgroups
$F(G)$ and $H$ are nilpotent. In view of $F(G)\trianglelefteq\ G$ and
Corollary~2.6.2 we conclude that $F(G)$ and $H$ are strongly permuteral in $G$.
 
$(2)\Rightarrow (3)$ obviously. 

$(3)\Rightarrow (1)$.
Let $\frak F =(G\ |\ $ the group $G=AB$ is the product of permuteral nilpotent
subgroups $A$ and $B$ of $G$).

Let $G$ be a group of minimal order in
$\frak F\setminus \frak U$. By the Theorem of Kegel-Wilandt \cite{Keg}, \cite{Wiel}
$G$ is soluble.
Let $N$ be a minimal normal subgroup of $G$. Since
$AN/N\simeq A/A\cap N\in\frak N$, $BN/N\simeq B/B\cap N\in\frak N$,
by Lemma~2.2(1) and by the choice of $G$ it follows that $G/N\in\frak U$.
Therefore $N$ is an unique minimal
normal subgroup of $G$ and $\Phi(G)=1$.
There exists a maximal subgroup $M$ in $G$ such that $G=NM$,
$N\cap M=1, {\rm {Core}}_G(M)=1$. $N=C_G(N)$ and $N$ is an elementary abelian $p$-group
for some $p\in\pi(G)$.
By Lemma~1.13, if $N\leq A$ then $A$ is a Sylow $p$-group
and $B$ is a Hall $p'$-group of $G$.

Consider any $x\in G$ for which $\langle x\rangle
B=B\langle x\rangle$ and $x\notin B$. Let $R=\langle x\rangle B$
and $\langle x\rangle=
\langle z\rangle\langle y\rangle$, where $\langle z\rangle\in{\rm{Syl}}_p
(\langle x\rangle)$ and $\langle y\rangle$ is a Hall $p'$-subgroup of
$\langle x\rangle$.
Clearly, $\langle x\rangle$ is not a
$p'$-group. By Lemma~1.14 $\langle y\rangle B$ is a Hall
$p'$-subgroup of $R$. So $\langle y\rangle\leq B$
and $R=\langle z\rangle B$.

Let
$B=P_1 \cdots P_k$, where $P_i\in
{\rm{Syl}}_{p_i}(B), i=1,  \ldots, k$. Take any
$i\in\{1,  \ldots, k\}$. By Lemma~1.14
for
$\pi_i=\{p, p_i\}$ there exists a Hall
$\pi_i$-subgroup $\langle z\rangle P_i=P_i\langle z\rangle$ in $R$.
Then $\langle x\rangle P_i=\langle z\rangle\langle y\rangle P_i=
\langle z\rangle P_i\langle y\rangle=P_i\langle z\rangle\langle y\rangle=
P_i\langle x\rangle$. Hence $x\in P_G(P_i)$. Since
$B\leq N_G(P_i)\leq P_G(P_i)$, we have 
$G=P_G(B)\leq P_G(P_i)$, i.e. $G=P_G(P_i)$. Since
$A\in {\rm{Syl}}_p(G), P_i\in {\rm{Syl}}_{p_i}(G)$, in view of Lemma~2.1(2) we conclude
that $G=P_G(H)$ for every Sylow subgroup $H$
of $G$.

By Lemma~2.5 $G$ is Ore dispersive. Since $N$ is an unique minimal normal subgroup of $G$
and $N$ is a
$p$-group, we obtain $p$ is the largest prime divisor
$|G|$. Then since $N\leq A\in {\rm{Syl}}_p(G)$, it follows that
$A \trianglelefteq\ G$. By Lemma~1.1 $G/C_G(N)=G/N\simeq M$
has not nonidentity normal $p$-subgroups. Therefore
$A\cap M=1$ and $N=A$. Since $G=AB=NB$, it follows that $B$ is a maximal subgroup of $G$.
In view of $P_G(B)=G$ and
$B\not= G$, there exists
$g\in G$ such that $\langle g\rangle
B=B\langle g\rangle$ and $g\notin B$ . Then $G= \langle g\rangle B$.
So $A$ is cyclic and $|A|=p$.
Now then $G/C_G(A)=G/A$ can be embedded in 
${\rm {Aut}}(Z_{p})\simeq Z_{p-1}$. Hence $G$ is supersoluble,
a contradiction. 
\hspace{\stretch{1}}$\square$

\medskip

{\bf Corollary 3.4.1.} {\it Let $G$ be a group, and let $G=AB$ be a product of its Sylow subgroups
$A$ and $B$. Then $G$ is supersolvable if and only if
$A$ and $B$ is permuteral in $G$. }

\medskip

{\bf Corollary 3.4.2.} {\it Let $G$ be a group. Then  $G$ is supersoluble if and only if
$G=F(G)H$, where $H$ is a permuteral Carter subgroup of $G$.}

\renewcommand{\refname}{\hfill\normalsize\bf References\hfill}

A.\,F.\,Vasil'ev, V.\,A.\,Vasil'ev

Francisk Skorina Gomel State University,
Sovetskaya str., 104,
Gomel 246019, Belarus.
E-mail address: formation56amail.ru, VovichX@mail.ru

\medskip

Т.\,I.\,Vasil'eva

Belarusian State University of Transport,
Kirov str., 34,
Gomel 246653, Belarus.
E-mail address: tivasilyeva@mail.ru


\end{document}